\documentclass{article}

\usepackage{amsfonts}
\usepackage{amssymb}

\usepackage[usenames,dvipsnames]{color}
\usepackage{amsmath}
\usepackage{amsfonts}
\usepackage{amssymb}
\usepackage{algorithmic,algorithm}
\usepackage{eurosym}
\usepackage{calc}
\usepackage{hyperref}
\usepackage{booktabs}

\usepackage{graphicx}
\usepackage{caption}
\usepackage{subcaption}
\usepackage{placeins}

\usepackage[margin=1in]{geometry}

\usepackage{amsthm}

\newtheorem{lemma}{Lemma}

\newtheorem{example}{Example}
\newtheorem{definition}{Definition}
\newtheorem{assumption}{Assumption}

\newtheorem{proposition}{Proposition}

\newcommand{\R}{\mathbb{R}}
\newcommand{\X}{\mathcal{X}}

\newcommand{\G}{G}
\newcommand{\E}{\mathcal{E}}
\newcommand{\V}{\mathcal{V}}

\newcommand{\I}{\mathbf{I}_n}

\title{On monotonicity of FIFO-diverging junctions}

\author{Marius Schmitt, John Lygeros}

\begin{document}

\maketitle

\begin{abstract}
This technical note concerns the dynamics of FIFO-diverging junctions in compartmental models for traffic networks. Many strong results on the dynamical behavior of such traffic networks rely on monotonicity of the underlying dynamics. In road traffic modeling, a common model for diverging junctions is based on the First-in, first-out principle. These type of junctions pose a problem in the analysis of traffic dynamics, since their dynamics are not monotone with respect to the positive orthant. However, this technical note demonstrates that they are in fact monotone with respect to the partial order induced by a particular, polyhedral cone.
\end{abstract}

\section{Introduction} \label{sec:introduction} 

Monotone systems are systems that preserve the ordering of trajectories, in particular, systems monotone with respect to the order induced by the positive orthant preserve the component-wise ordering of trajectories \cite{hirsch2006monotone}. Monotonicity of the dynamics of compartmental systems, with respect to the positive orthant, has been used widely in analyzing the dynamical behavior of certain traffic networks \cite{gomes2008behavior,lovisari2014stability2,coogan2016stability,como2017resilient}. In particular, it has been shown that monotone routing policies show favorable resilience to capacity reductions \cite{como2013robust1,como2013robust2} and that such policies can be used to stabilize maximal-throughput equilibria \cite{como2015throughput}. However, it is well known that the dynamics of First-in, first-out (FIFO) diverging junctions are \emph{not} monotone with respect to the positive orthant, since congestion in any downstream cell can block flows into other downstream cells \cite{munoz2002bottleneck,coogan2014dynamical}. Different, monotone models for diverging junctions have been suggested \cite{lovisari2014stability2}. However, these models do \emph{not} preserve the turning rates and hence, they are not suitable for the Freeway Network Control (FNC) problem, where turning rates are assumed to be constant. In addition, there is strong empirical evidence for FIFO-behavior of diverging junctions \cite{munoz2002bottleneck}. It has been shown that the dynamics of FIFO diverging junctions satisfy a mixed-monotonicity property \cite{coogan2016stability}, that is, they can be embedded into a higher-dimensional monotone system. However, this property is somewhat weaker than monotonicity. 

Alternatively, one can ask whether the dynamics of FIFO-diverging junctions are monotone with respect to a different, partial order. It is known they are not monotone with respect to any partial order induced by an orthant \cite{coogan2016stability}. The main purpose of this technical note is to show that the dynamics of FIFO diverging junctions are monotone with respect a polyhedral cone, defined in the following, that is not an orthant. \\

We use the following notation: the symbols $\geq, \leq$ denote component-wise inequalities. Generalized inequalities with respect to some closed, convex and pointed cone $K$ are denoted by $\succeq_K$ and $\preceq_K$. The closed, positive orthant is denoted $\R^n_+$. Other sets will be denoted using calligraphic letters, e.g.\ $\mathcal{V}$. Notation for describing the compartmental traffic model, and the associated graph on which it is defined, is introduced in the following section.

\section{System model} \label{sec:model} 
First-order compartmental models are widely used to model the evolution of traffic networks \cite{daganzo1994cell,daganzo1995cell,gomes2006optimal,gomes2008behavior,coogan2016stability}. Here, we consider a compartmental model based on a directed graph $\G = (\V, \E)$ with $\E \subset \V \times \V$. The vertices $v \in \V$ model junctions, while the edges $e \in \E$ model cells between junctions. This technical note focusses on the FIFO model for diverging junctions and therefore, we restrict our attention to graphs $\G$, which are rooted, directed trees. Such a tree has a unique root $v_r \in \V$ with $\deg^-( v_r ) = 0$, where $\deg^-(v)$ denotes the in-degree of a vertex. All other junctions $v \neq v_r$ have $\deg^-(v) = 1$. We assume that $\deg^+( v_r ) = 1$, where $\deg^+(v)$ denotes the out-degree of a vertex, but allow an arbitrary out-degree for all other vertices. The unique cell originating at the root (-vertex) $v_r$ is denoted by $r \in \E$. Vertices with $\deg^+(v) = 0$ are called sinks. The head of edge (cell) $e$ is denoted by $\sigma_e$ and the tail by $\tau_e$. Traffic flows from tail $\tau_e$ to head $\sigma_e$. An example network is depicted in Figure \ref{fig:ex_topology}.

The state of the compartmental model is comprised of the states $x_e(t) \in [0, \bar x_e]$ of the individual cells, where $\bar x_e$ is the jam density of edge $e$. To describe the evolution of traffic, we assume that every edge is equipped with a demand function, modeling the amount of traffic that seeks to travel downstream and a supply function, modeling the available, free space.
\begin{assumption} \label{assumption:fd} 
Every demand function $d_e( x_e ) : [0, \bar x_e] \to \R^n_+$ is nondecreasing, Lipschitz continuous and $d_e ( 0 ) = 0$. Every supply function $s_e( x_e ) : [0, \bar x_e] \to \R^n_+$ is nonincreasing, Lipschitz continuous and $s_e( \bar x_e ) = 0$. 
\end{assumption}
\begin{figure}[t] 
	\centering
	\includegraphics[width=6cm]{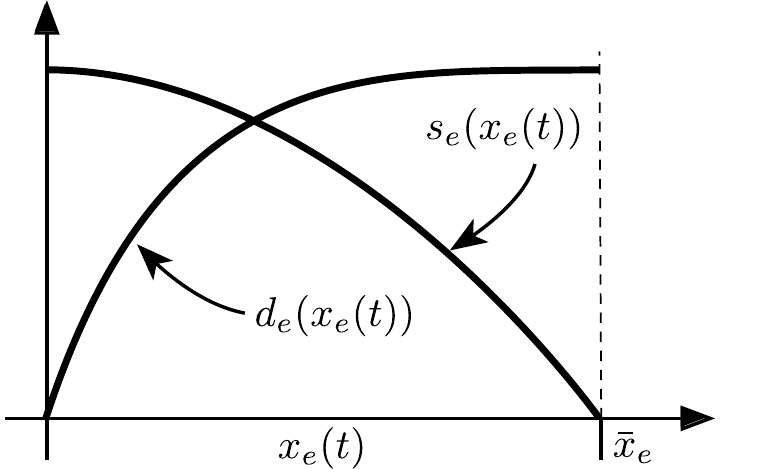}
	\caption{Demand and supply function satisfying Assumption \ref{assumption:fd}.}
	\label{fig:fd}
\end{figure} 

Examples of demand and supply functions are depicted in Figure \ref{fig:fd}. Demand and supply functions are used to model maximal cell outflow and inflow, respectively. For junctions with $\deg^-(v) > 1$, one needs to define how flow leaving an upstream cell is distributed onto downstream cells. The percentage of flow leaving cell $i$ that is routed to cell $e$ is described by the constant turning rates $\beta_{e,i}$,  for any pair of adjacent cells $e$, $i$ where $\sigma_i = \tau_e$. Conservation of traffic requires that $\sum_{i \in \E} \beta_{e,i} \leq 1$. In case $\sum_{i \in \E} \beta_{e,i} < 1$, we assume that the remaining flow has left the modeled part of the traffic network. The FIFO model for diverging junctions, that is, for junctions with $\deg^+(v) \geq 2$, maximizes traffic flows such that the total flow leaving a cell is bounded by its demand function and the flow entering any cell is bounded by its supply function. With $\phi_e(t)$ denoting the flow out of cell $e$ and $\phi_e^{in}(t)$ the flow into this cell, the evolution of the compartmental model is given as
\begin{equation}
\label{eq:system}
\dot x(t) = f \big( t, x(t) \big) := \phi^{in}(t) - \phi(t) ,
\end{equation}
with
\begin{equation}
\label{eq:flow}
\phi_e(t) = \min \Bigg\{ d_e \big( x_e(t) \big),~ \min_{i: \beta_{i,e} > 0} \frac{ s_i \big( x_i(t) \big)}{\beta_{i,e} } \Bigg\}
\end{equation}
and 
\begin{equation}
\label{eq:inflow}
\phi_e^{in}(t) = \begin{cases} \beta_{e,i} \phi_i(t), & \forall e \neq r, ~ \beta_{e,i} > 0, \\ \min \Big\{ w_r(t), s_e \big( x_e(t) \big) \Big\}, & e = r, \end{cases}
\end{equation}
Note that in \eqref{eq:inflow}, there exists a unique upstream cell $i$ for every cell $e$, such that $\beta_{e,i} > 0$, since the network graph is a rooted tree. Here, the quantity $w_r(t)$ denotes external demand for the cell incident to the source. If this external inflow exceeds supply of the cell, surplus external demand is discarded.\footnote{An alternative model for source cells ensures that all external inflow is served by postulating that source cells have infinite capacity.} Equation \eqref{eq:flow} encodes the FIFO property: the flow $\phi_e$ is limited by the minimum of the scaled supply among all downstream cells. Depletion of free space in any downstream cell also limits flow from the upstream cell in all other downstream cells. The solution to the system described by \eqref{eq:system}, \eqref{eq:flow} and \eqref{eq:inflow} is well defined for all $t \geq 0$ and the convex set $\X = \prod_{e \in \E} [0, \bar x_e]$ is forward-invariant as the system is a special case of the system described in \cite{coogan2016stability}. In the following, we denote the solution of the compartmental model as $\Psi_f(t,x_0)$. 

\section{Monotonicity} \label{sec:monotonicity} 

Cones are instrumental in defining partial orders, and in turn, monotone systems. For completeness, we state the basic definition of (pointed, proper) cones (according to \cite{boyd2004convex}).

\begin{definition}[Cone] 
A set $K \subseteq \R^n$ is a \emph{cone} if for every $x \in K$ and every $\lambda \in \R$, $ \lambda > 0$, it follows that $\lambda x \in K$. A cone is called \emph{pointed} if $K \cap (-K) = \{ 0 \}$. It is called \emph{proper} if it is closed, convex, pointed and has non-empty interior.
\end{definition} 

Closed, convex and pointed cones $K \subset \R^n$ are of particular interest, since such a cone can be used to define a partial order on $\R^n$, via generalized inequalities \cite[Proposition 3.38 ]{rockafellar2009variational}. That is, the cone $K$ induces a partial order, where $x \succeq_K y$ iff $x - y \in K$.\footnote{Some authors, in particular \cite{boyd2004convex}, restrict their attention to proper cones, that is, closed, convex and pointed cones \emph{with non-empty interior}, when defining a partial order. Assuming non-empty interior is not necessary for defining a partial order itself \cite[Proposition 3.38 ]{rockafellar2009variational}, but this assumption is made in infinitesimal characterizations of monotone systems, such as the quasi-monotone conditions \cite[Theorem 3.2]{hirsch2006monotone}..} Such a partial order is used to define a monotone system \cite{hirsch2006monotone}.

\begin{definition} \label{definition:monotone} 
An autonomous system $\Psi : \R_+ \times \X \to \X$ is \emph{monotone} with respect to the (closed, convex and pointed) cone $K$ if for all $x_0, y_0 \in \X$ with $x_0 \succeq_K y_0$ it holds that
\begin{equation*}
\Psi(t,x_0) \succeq_K \Psi(t,y_0)
\end{equation*}
for all $t \geq 0$.
\end{definition} 

If the cone equals the closed, positive orthant $K = \R^n_+$, the standard, componentwise inequality $x \succeq_{\R^n_+} y ~\Leftrightarrow~ x \geq y$ is obtained. There also exist infinitesimal characterizations of monotonicity such as the \emph{quasi-monotone condition} \cite[Theorem 3.2]{hirsch2006monotone}. For systems that are monotone with respect to the positive orthant, the quasi-monotone condition reduces to the \emph{Kamke-M\"uller conditions} \cite[Equation 3.3]{hirsch2006monotone}.

\begin{lemma}[Kamke-M\"uller conditions] 
An autonomous system of the form \eqref{eq:system} is monotone with respect to $\R^n_+$ iff
\begin{equation*}
x \leq y \text{ and } x_i = y_i ~\implies~ f_i(t,x) \leq f_i(t,y)
\end{equation*}
\end{lemma} 

The Kamke-M\"uller condition means that $f_i(t,x)$ is nondecreasing in $x_j$ for $i \neq j$. The dynamics of a freeway segment with only onramp and offramp junctions are monotone with respect to the positive orthant, a property that can be leveraged to analyze its stability properties \cite{gomes2008behavior}. In addition, widely-used merging junction models exhibit monotone (w.r.t.\ $\R^n_+$) dynamics \cite{lovisari2014stability2,coogan2014dynamical,coogan2016stability}. However, the dynamics of FIFO-diverging junctions are \emph{not} monotone with respect to the positive orthant, since congestion in any one downstream cell can block flows into other downstream cells, thereby decreasing the density in those cells \cite{munoz2002bottleneck,coogan2014dynamical,coogan2016stability}. It has also been proven that the dynamics of FIFO diverging junctions are not monotone with respect to any orthant \cite{coogan2016stability}, a result based on the ``graphical condition" according to \cite[Proposition 2]{angeli2004interconnections}.

The main purpose of this technical note is to show that the dynamics of the network introduced in Section \ref{sec:model}, and hence the dynamics of FIFO diverging junctions, are monotone with respect to the order induced by a polyhedral cone, defined in the following, that is not an orthant. To do so, consider the routing matrix $R$ with entries $R_{i,e} := \beta_{i,e}$, whenever the turning rate is defined, and $R_{i,e} = 0$ otherwise. The routing matrix is column-substochastic, that is, its column sums are smaller than or equal to one, $\sum_{i \in \E} R_{i,e} \leq 1$ for all $e \in \E$, because of the conservation law of traffic. In addition, its spectral radius is strictly smaller than one, $\rho( R ) < 1$, for directed-tree networks as considered in this note. The latter fact is equivalent to the assertion that all traffic eventually leaves the network \cite{varaiya2013max}. Using the routing matrix, we can write the system dynamics as
\begin{align*}
\dot x = (R - \I) \phi(t) + \hat e_r \phi^{in}_r(t) ,
\end{align*}
where $\hat e_r$ is the unit vector for which the component corresponding to cell $r$ is equal to one. Consider now $P = (\I - R)^{-1} = \sum_{k=0}^{\infty} R^k$, which has nonnegative entries $P \geq 0$ since $R \geq 0$ \cite{coogan2016stability}. 

\begin{proposition} \label{proposition:fifo} 
The dynamics of the compartmental model \eqref{eq:system}-\eqref{eq:inflow} satisfying Assumption \ref{assumption:fd}, defined on a directed tree with FIFO-diverging junctions, are monotone with respect to the partial order induced by the polyhedral cone $\mathcal{P} := \{x : Px \geq 0\}$.
\end{proposition} 

For now, we will avoid employing the quasi-monotone condition and instead performe a state transformation and prove monotonicity, with respect to $\R^n_+$, of the transformed system.
\begin{proof} 
Consider the state transformation $z(t) := P x(t)$ and the transformed system
\begin{align*}
\dot z = g \big(t, z(t) \big) := P \cdot f \big( t, P^{-1} z(t) \big) = (\I - R)^{-1} \cdot \big( (R - \I) \phi(t) + \hat e_r \phi^{in}_r(t) \big) = -\phi(t) + \hat e_r \phi^{in}_r(t) ,
\end{align*}
In the last equality, we have used that $P \hat e_r = P (\I - R) \hat e_r = \hat e_r$. The transformed system is defined on the convex set $\mathcal{Z} := \{ z : (\I-R) z \in \mathcal{X} \}$. We have that $x_e(t) = z_e(t) - \beta_{e,e^+} z_{e^+}(t)$, where $e^+$ is the unique cell upstream of $e$, for all $e : \tau_e \neq r$, and $x_e(t) = z_e(t)$ for $e : \tau_e = r$. The transformed system is monotone with respect to the positive orthant, which can be verified via the Kamke-M\"uller conditions. Recall that demand and supply functions are Lipschitz-continuous. From the equations defining the compartmental model, it follows that $f(t,x)$ according to \eqref{eq:system}, and in turn $g(t,z)$, are Lipschitz-continuous. Lipschitz continuity implies that the components $g_i(t,z)$ are differentiable almost everywhere, and hence, verifying that 
\begin{equation*}
\frac{ \partial g_i }{ \partial z_j } (z) \geq 0 \quad \forall i \neq j, ~ \forall z \in \mathcal{Z} ,
\end{equation*}
whenever the partial derivative exists, is sufficient for the Kamke-M\"uller conditions to hold. In the following, we will take partial derivatives of expressions involving demand and supply functions, whereby we implicitly assume that the corresponding expressions hold, whenever the partial derivative exists. We first establish that 
\begin{align*}
\frac{ \partial }{\partial z_j} \phi_r^{in}(t) = \frac{ \partial }{\partial z_j} \cdot \min \big\{ \phi_r(t), s_r( z_r ) \big\} \geq \begin{cases} s_r' ( z_r ), & j = r, \\ 0, & \text{else}, \end{cases}
\end{align*}
almost everywhere in $\mathcal{Z}$. Therefore, for $j \neq e$,
\begin{align*}
\frac{\partial }{\partial z_j} g_e\big(t,  z(t) \big) 
 &\geq  -\frac{\partial}{\partial z_j} \min \Big\{ d_e \big( x_e(z) \big),~ \min_{i: \beta_{i,e} > 0} \beta_{i,e}^{-1} s_i \big( x_i(z) \big) \Big\} + \frac{ \partial }{\partial z_j} \phi_r^{in}(t) \\
 &\geq - \max \bigg\{ \frac{\partial}{\partial z_j}  d_e \big( x_e(z) \big),~ \max_{i: \beta_{i,e} > 0} \beta_{i,e}^{-1} \frac{\partial}{\partial z_j}  s_i \big( x_i(z) \big) \bigg\} \\
 &\geq \min \bigg\{ -\frac{\partial}{\partial z_j} d_e \big( x_e(z) \big),~ \min_{i: \beta_{i,e} > 0} -\beta_{i,e}^{-1} \cdot \frac{\partial}{\partial z_j} s_i \big( z_i - \beta_{i,e} z_{e} \big) \bigg\} ,
\end{align*}
almost everywhere. We consider the partial derivatives individually and find that for $j \neq e$,
\begin{align*}
- \frac{\partial}{\partial z_j} d_e \big( x_e(z) \big) = \begin{cases} \beta_{e,e^+} \cdot d_e'(z_e - \beta_{e,e^+} z_{e^+}) \geq 0, & j = e^+, \\ 0, & \text{else}, \end{cases}
\end{align*}
where we have used that the demand function is nondecreasing. Furthermore, for $j \neq e$,
\begin{align*}
-\beta_{i,e}^{-1} \cdot \frac{\partial}{\partial z_j} s_i \big( z_i - \beta_{i,e} z_{e} \big) = \begin{cases} -\beta_{i,e}^{-1} \cdot s_i' ( z_i - \beta_{i,e} z_{e} ) \geq 0, & j = i ,\\ 0, & \text{else}, \end{cases}
\end{align*}
where we have used that the supply function is nonincreasing. Hence,
\begin{align*}
\frac{\partial g_e\big( z \big)}{\partial z_j} \geq 0, \quad\quad & \forall j \neq e,
\end{align*}
almost everywhere in $\mathcal{Z}$, which implies that the transformed system is monotone with respect to the positive orthant. Monotonicity of the transformed system means that for all $z_0, w_0 \in \mathcal{Z}$, the implication $z_0 \geq w_0 \implies \Psi_g(t,z_0) \geq \Psi_g(t,w_0)$ holds true for all $t \geq 0$. In turn, this means that for all $x_0, y_0 \in \X$,
\begin{align*}
x_0 \succeq_{K} y_0 ~\Leftrightarrow~ P ( x_0 - y_0 ) \geq 0  &~\implies~ \Psi_g \big( t, P x_0 \big) \geq \Psi_g \big( t, P y_0 \big)
\end{align*}
and
\begin{align*}
\Psi_g \big( t, P x_0 \big) \geq \Psi_g \big( t, P y_0 \big) &\Leftrightarrow~ P \cdot \Psi_f \big( t, x_0 \big) \geq P \cdot \Psi_f \big( t, y_0 \big) \\
~&\Leftrightarrow~ \Psi_f(t,x_0) \succeq_{\mathcal{P}} \Psi_f(t,y_0) ,
\end{align*} 
which proves monotonicity of the original compartmental model with respect to $\mathcal{P} = \{x : Px \geq 0\}$.
\end{proof} 

\section{Discussion} \label{sec:discussion} 

One can avoid to introduce a state transformation and verify the quasi-monotone condition directly. Note that the cone $K$ is proper, which is required for applying the quasi-monotone conditions according to \cite[Theorem 3.2]{hirsch2006monotone}. In this case, one needs to verify that for all $x,y \in \X$, 
\begin{align*}
x \succeq_{\mathcal{P}} y ,~ \zeta^\top x = \zeta^\top y ~\implies~ \zeta^\top f(x) \geq \zeta^\top f(y)
\end{align*}
for all $\zeta \in \mathcal{P}^*$, where $\mathcal{P}^* = \big\{ \zeta : \zeta = P^\top \lambda,~ \lambda \geq 0 \big\} $ is the dual cone.\footnote{Strictly speaking, the dual cone is a subset of the linear functions that operate on elements of $\mathcal{P}$, but its elements $x \mapsto \zeta^\top x$ can be identified with $\zeta$.} Since every element of the dual cone can be obtained as a positive, linear combination of the row vectors $P_{(e,:)}$ of $P$, this condition is equivalent to verifying that $x \succeq_{\mathcal{P}} y$ and $P_{(e,:)} x \geq P_{(e,:)} y$ implies that $P_{(e,:)} f(x) \geq P_{(e,:)} f(y)$, which reduces to verifying the Kamke-M\"uller conditions of the transformed system.

One advantage of introducing the transformed system is that the states $z_e(t)$ admit an intuitive explanation. Consider the case when no further, external flows enter the network from time $t$ onwards, that is, $\phi^{in}_r(\tau) = 0$ for $\tau \geq t$. Then,
\begin{align*}
z_e(t) = z_e(+\infty) + \int_{+\infty}^{t} \dot z~ \mathrm{d}\tau = \int^{+\infty}_{t} \phi_e(\tau) ~\mathrm{d}\tau 
\end{align*}
where $z( +\infty ) = P x( +\infty ) = 0$, because all traffic eventually leaves the network, by virtue of the network structure encoded in $R$. Hence, one can interpret $z_e(t)$ as the \emph{cumulative traffic demand} that has to be served by cell $e$ in the future, assuming no further, external traffic enters the network. For a network graph that is a directed tree as assumed in this technical note, $z_e$ is simply the sum of $x_e(t)$ and all traffic volume $x_i$ in cells $i$ upstream of cell $e$, weighted according to what percentage of upstream traffic volume will eventually be routed to cell $e$.


\begin{figure}[t] 
	\centering
		\begin{subfigure}[b]{0.58\textwidth}
		\includegraphics[width=\textwidth]{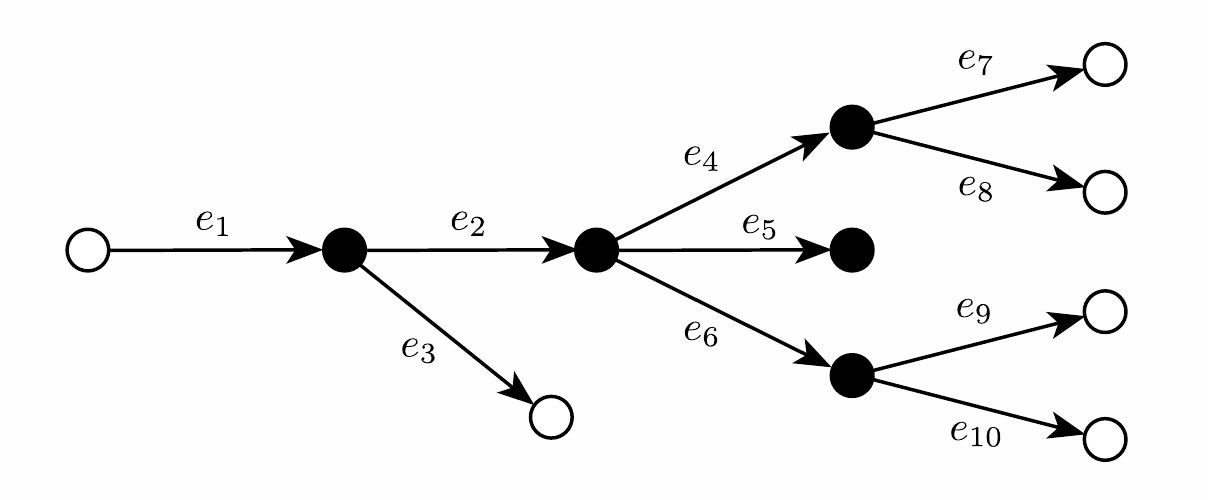}
		\caption{Network topology.}
		\label{fig:ex_topology}
	\end{subfigure} \hspace{1cm}	
	\begin{subfigure}[b]{0.31\textwidth}
\begin{center}
\begin{tabular}{rl}
\toprule 
$k$ & Cell indices \\ \hline
1 & $1:10$ \\
2 & $\{ 1,2,3,4,10 \}$ \\
3 & $\{ 1,6,7 \}$ \\
4 & $\{ 2,7,9,10 \}$ \\
5 & $\{ 4,5,6 \}$ \\
6 & $\varnothing$ \\
\bottomrule
\end{tabular}
\end{center}
		\caption{Indices of initially congested cells}
		\label{fig:ex_table}
	\end{subfigure} \\[1ex]
	\begin{subfigure}[b]{0.31\textwidth}
		\includegraphics[width=\textwidth]{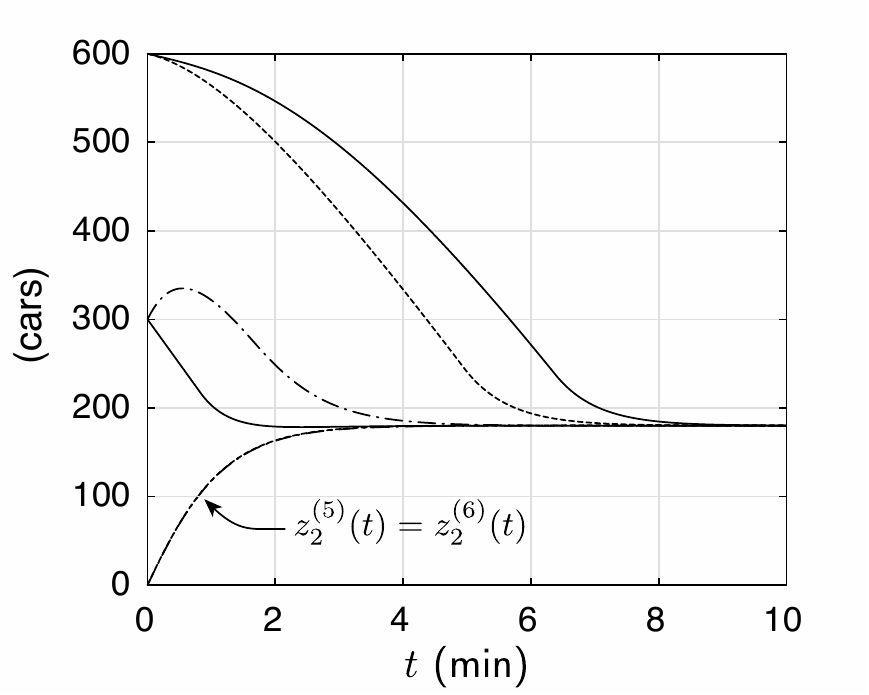}
		\caption{Evolution of $z_2^{(k)}(t)$.}
		\label{fig:ex_cell2}
	\end{subfigure} ~	
	\begin{subfigure}[b]{0.31\textwidth}
		\includegraphics[width=\textwidth]{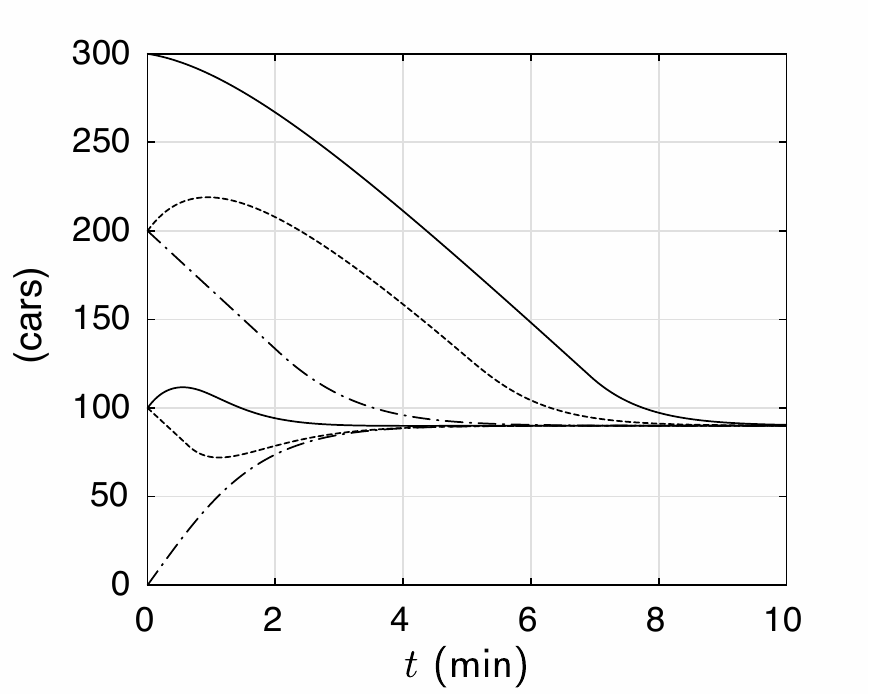}
		\caption{Evolution of $z_6^{(k)}(t)$.}
		\label{fig:ex_cell6}
	\end{subfigure} ~
	\begin{subfigure}[b]{0.31\textwidth}
		\includegraphics[width=\textwidth]{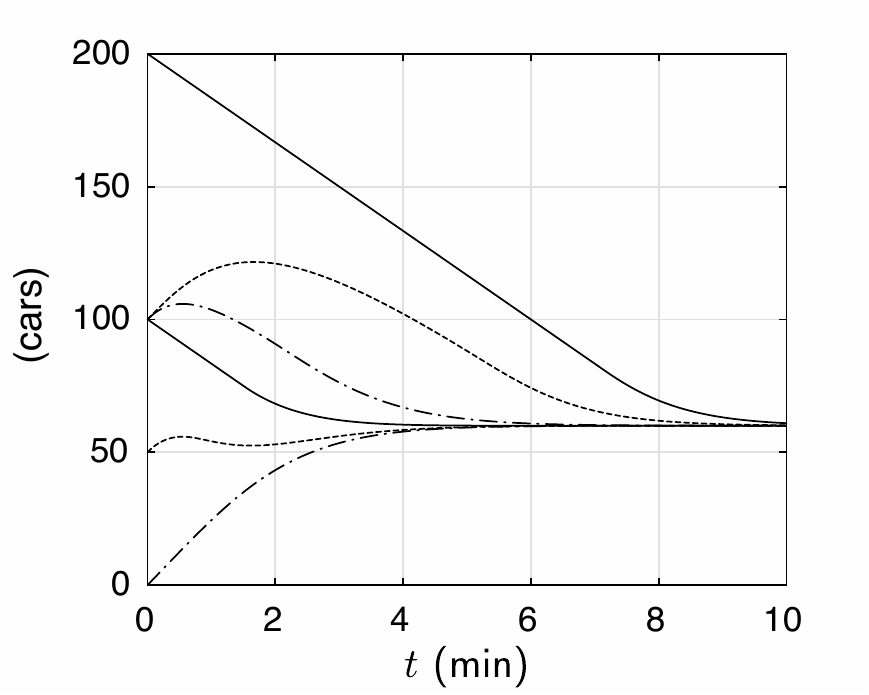}
		\caption{Evolution of $z_9^{(k)}(t)$.}
		\label{fig:ex_cell9}
	\end{subfigure} 
	\caption{The network topology (a) is used for simulations, where different cells are congested in the initial state according to (b).
	Figures (c)-(e) depict the evolution of cumulative, future flows $z_e^{(k)}(t)$ for certain cells. Note that the trajectories retain their ordering.    }
	\label{fig:example}
\end{figure} 

To verify Proposition \ref{proposition:fifo} numerically, we also present the following, simple example.
\begin{example} \label{ex:network} 
We simulate the example network depicted in Figure \ref{fig:ex_topology} for different initial conditions $x^{(k)}(0)$. The turning rates are $\beta_{2,1} = 0.9$, $\beta_{3,1} = 0.1$, $\beta_{4,2} = \beta_{5,2} = \beta_{6,2} = 1/3$ and $\beta_{7,4} = \beta_{8,4}  = \beta_{9,6} = \beta_{10,6} = 1/2$. Demand $d_e( x_e ) = \min \{ v_e x_e, F_e \}$ and supply functions $s_e( x_e ) = \min \{ F_e, w_e ( x_e - \bar x_e ) \}$ are piecewise-affine, with $v_e = 100$ and $w_e = 100/3$. The cell capacities $F_e$ are chosen such that all cells reach their capacity limit simultaneously, for steady-state flows with $\phi^{in}_1 = 50000/3$.\footnote{Such a choice ensures that any congested cell can obstruct upstream demand, such that the FIFO-diverging dynamics come into effect, but the results do not depend on such a choice. Similar results are obtained if the cell capacities are randomly disturbed, and the initial densities are adapted accordingly, such that they are still ordered with respect to $\mathcal{P}$.} Note that for this example, all parameters and quantities are dimensionless. In the initial states, certain cells are congested $x_e^{(k)}(0) = 2 \cdot \frac{F_e}{v_e}$, while the remaining cells are empty $x_e^{(k)}(0) = 0$. The cells that are congested for each $k$ are listed in Figure \ref{fig:ex_table}. It can be verified that the initial states are ordered in the sense that $x^{(1)}(0) \succeq_{\mathcal{P}} x^{(2)}(0) \succeq_{\mathcal{P}} \dots \succeq_{\mathcal{P}} x^{(6)}(0) $. Monotonicity implies that the ordering of trajectories is preserved, which can be visually verified by depicting the transformed states $z_e^{(k)}(t)$ for different cells and confirming that the trajectories do not intersect. This is indeed observed in Figures \ref{fig:ex_cell2}-\ref{fig:ex_cell9}. 
\end{example} 

A natural follow-up question to the result in this note is to ask whether the dynamics of typical models for merging junctions, many of which are monotone with respect to the positive orthant, are also monotone with respect to the order induced by $\mathcal{P}$. Unfortunately, it can be shown by counterexample that this is \emph{not} the case for the most important merging models, Daganzo's priority rule \cite{daganzo1995cell} and the proportional-priority merging model \cite{kurzhanskiy2010active,coogan2016stability}. This means that compartmental models for traffic networks, that contain both FIFO-diverging junctions and merging junctions of one of the two described types are neither monotone with respect to $\R^n_+$ nor $\mathcal{P}$, which limits the applicability of monotone system theory in the study of traffic networks. However, a partial remedy is described in \cite{schmitt2017exact}: if merging flows are controlled, one can recover monotonicity of the dynamics of such a merging junction.\footnote{Strictly speaking, \cite{schmitt2017exact} uses a reformulation of the system dynamics related to the state transformation used in the proof of Proposition \ref{proposition:fifo}, which turns out to be monotone with respect to the positive orthant if the control inputs are held constant and the resulting, autonomous system is considered.} In fact, it turns out that monotonicity is crucial in deriving a tight, convex relaxation of the FNC problem for the corresponding traffic network.

\section{Conclusions} \label{sec:conclusions} 
In this technical note, we have demonstrated that FIFO-diverging junctions are monotone with respect to the partial order induced by $\mathcal{P}$. Furthermore, we have seen that this ordering can be interpreted as being based on the cumulative, future traffic flow. However, typical dynamics of merging junctions, which are known to be monotone with respect to $\R^n_+$, are not monotone with respect to the ordering induced by $\mathcal{P}$. This means that while tools from monotone system theory can in principle be applied to compartmental models with only FIFO-diverging junctions, networks that contain \emph{both} FIFO-diverging junctions and merging junctions described Daganzo's priority rule or the proportional-priority merging model still present a challenge. This is, of course, a major limitation. So far, the main application of the results in this note and, in fact, our main motivation in pursuing this research, is their immediate application in the FNC problem with \emph{controlled merging flows}, where additional assumptions on the available actuation help to restore monotonicity of the dynamics of merging junctions.

\FloatBarrier

\bibliography{/Users/mschmitt/Documents/Docs_Latex/MS_bibliography}

%

\end{document}